\documentclass[a4paper,pdftex,reqno,10pt]{amsart}

\allowdisplaybreaks[1]

\usepackage{geometry}

\usepackage{graphicx}
\usepackage{enumerate}
\usepackage{courier}
\usepackage{times}
\usepackage{amsmath,amssymb}

\usepackage{hyperref}

\theoremstyle{plain}
\newtheorem{Theorem}{Theorem}[section]
\newtheorem{Proposition}[Theorem]{Proposition}

\newtheorem{Lemma}[Theorem]{Lemma}

\newenvironment{Proof}
{\begin{trivlist}\item[]{{\sc Proof.}}}{\hfill{$\square$}\noindent\end{trivlist}}

\theoremstyle{definition}
\newtheorem{Definition}[Theorem]{Definition}

\theoremstyle{remark}

\begin{document}


\title{The ${[46,9,20]_2}$ code is unique}


 \author{Sascha Kurz}
 \address{Sascha Kurz, University of Bayreuth, 95440 Bayreuth, Germany}
 \email{sascha.kurz@uni-bayreuth.de}

\abstract{The minimum distance of all binary linear codes with dimension at most eight is known. The smallest open case for 
dimension nine is length $n=46$ with known bounds $19\le d\le 20$. Here we present a $[46,9,20]_2$ code and show its uniqueness. 
Interestingly enough, this unique optimal code is asymmetric, i.e., it has a trivial automorphism group. Additionally, we show the 
non-existence of $[47,10,20]_2$ and $[85,9,40]_2$ codes.
}}

\maketitle

\section{Introduction}
An $[n,k,d]_q$-code is a $q$-ary linear code with length $n$, dimension $k$, and minimum Hamming distance $d$. Here we will only consider 
binary codes, so that we also speak of $[n,k,d]$-codes. Let $n(k,d)$ be the smallest integer $n$ for which an $[n,k,d]$-code exists. Due to 
Griesmer \cite{griesmer1960bound} we have
\begin{equation}
  \label{eq_griesmer_bound}
  n(k,d)\ge g(k,d):=\sum_{i=0}^{k-1} \left\lceil\frac{d}{2^i}\right\rceil,
\end{equation} 
where $\lceil{x}\rceil$ denotes the smallest integer $\ge x$. As shown by Baumert and McEliece \cite{baumert1973note} for every fixed dimension $k$ 
there exists an integer $D(k)$ such that $n(k,d)=g(k,d)$ for all $d\ge D(k)$, i.e., the determination of $n(k,d)$ is a finite problem for every 
fixed dimension $k$. For $k\le 7$, the function $n(k,d)$ has been completely determined by Baumert and McEliece \cite{baumert1973note} and 
van Tilborg \cite{van1981smallest}. After a lot of work of different authors, the determination of $n(8,d)$ has been completed by 
Bouyukliev, Jaffe, and Vavrek \cite{bouyukhev2000smallest}. For results on $n(9,d)$ we refer e.g.\ to \cite{dodunekov1999some} and the references therein. 
The smallest open case for dimension nine is length $n=46$ with known bounds $19\le d\le 20$. Here we present a $[46,9,20]_2$ code and show its uniqueness. 
Interestingly enough, this unique optimal code is asymmetric, i.e., it has a trivial automorphism group. Speaking of a $\Delta$-divisible code for 
codes whose weights of codewords all are divisible by $\Delta$, we can state that the optimal code is $4$-divisible. $4$-divisible codes are also 
called doubly-even and $2$-divisible codes are called even. Additionally, we show the non-existence of $[47,10,20]_2$ and $[85,9,40]_2$ codes.

Our main tools -- described in the next section -- are the standard residual code argument (Proposition~\ref{prop_residual}), the MacWilliams identities 
(Proposition~\ref{prop_MacWilliams_identitites}), a result based on the weight distribution of Reed-Muller codes (Proposition~\ref{prop_div_one_more}), 
and the software packages \texttt{Q-Extension} \cite{bouyukliev2007q}, \texttt{LinCode} \cite{LinCode} to enumerate linear codes with a list of allowed 
weights. For an easy access to the known non-existence results for linear codes we have used the online database \cite{Grassl:codetables}. 

\section{Basic tools}

\begin{Definition}
  Let $C$ be an $[n,k,d]$-code and $c\in C$ be a codeword of weight $w$. The restriction to the support of $c$ is called the residual code 
  $\operatorname{Res}(C;c)$ of $C$ with respect to $c$. If only the weight $w$ is of importance, we will denote it by $\operatorname{Res}(C;w)$.
\end{Definition}

\begin{Proposition}
  \label{prop_residual}
  Let $C$ be an $[n,k,d]$-code. If $d>w/2$, then $\operatorname{Res}(C;w)$ has the parameters 
  $$
    \left[n-w,k-1,\ge d-\left\lfloor w/2\right\rfloor\right].
  $$    
\end{Proposition}
Some authors call the result for the special case $w=d$ the one-step Griesmer bound.

\begin{Proposition} (\cite{macwilliams1977theory}, MacWilliams Identities)
  \label{prop_MacWilliams_identitites}
  Let $C$ be an $[n,k,d]$-code and $C^\perp$ be the dual code of $C$. Let $A_i(C)$ and $B_i(C)$ be the number of codewords of weight $i$ in $C$ and 
  $C^\perp$, respectively. With this, we have
  \begin{equation}
    \label{eq_MacWilliams}
    \sum_{j=0}^n K_i(j)A_j(C)=2^kB_i(C),\quad 0\le i\le n
  \end{equation} 
  where
  $$
    K_i(j)=\sum_{s=0}^n (-1)^s {{n-j}\choose{i-s}}{{j}\choose{s}},\quad 0\le i\le n 
  $$
  are the binary Krawtchouk polynomials. We will simplify the notation to $A_i$ and $B_i$ whenever $C$ is clear from the context.
\end{Proposition}
Whenever we speak of the first $l$ MacWilliams identities, we mean Equation~(\ref{eq_MacWilliams}) for $0\le i\le l-1$. 
Adding the non-negativity constraints $A_i,B_i\ge 0$ we obtain a linear program where we can maximize or minimize certain quantities, which is called the 
linear programming method for linear codes. Adding additional equations or inequalities strengthens the formulation. 

\begin{Proposition}(\cite[Proposition 5]{dodunekov1999some}, cf.~\cite{simonis1994restrictions})
  \label{prop_div_one_more}
  Let $C$ be an $[n,k,d]$-code with all weights divisible by $\Delta:=2^a$ and let $\left(A_i\right)_{i=0,1,\dots,n}$ be the weight distribution of $C$. Put
  \begin{eqnarray*}
    \alpha&:=&\min\{k-a-1,a+1\},\\
    \beta&:=&\lfloor(k-a+1)/2\rfloor,\text{ and}\\ 
    \delta&:=&\min\{2\Delta i\,\mid\,A_{2\Delta i }\neq 0\wedge i>0\}.
  \end{eqnarray*}
  Then the integer 
  $$
    T:=\sum_{i=0}^{\lfloor n/(2\Delta)\rfloor} A_{2\Delta i}
  $$  
  satisfies the following conditions.
  \begin{enumerate}
    \item \label{div_one_more_case1}
          $T$ is divisible by $2^{\lfloor(k-1)/(a+1)\rfloor}$.
    \item \label{div_one_more_case2}
          If $T<2^{k-a}$, then
          $$
            T=2^{k-a}-2^{k-a-t}
          $$
          for some integer $t$ satisfying $1\le t\le \max\{\alpha,\beta\}$. Moreover, if $t>\beta$, then $C$ has an $[n,k-a-2,\delta]$-subcode and if $t\le \beta$, it has an 
          $[n,k-a-t,\delta]$-subcode.
    \item \label{div_one_more_case3}
          If $T>2^k-2^{k-a}$, then
          $$
            T=2^k-2^{k-a}+2^{k-a-t}
          $$        
          for some integer $t$ satisfying $0\le t\le \max\{\alpha,\beta\}$. Moreover, if $a=1$, then $C$ has an $[n,k-t,\delta]$-subcode. If $a>1$, then $C$ has an 
          $[n,k-1,\delta]$-subcode unless $t=a+1\le k-a-1$, in which case it has an $[n,k-2,\delta]$-subcode.
  \end{enumerate}
\end{Proposition}

A special and well-known subcase is that the number of even weight codewords in a $[n,k]$ code is either $2^{k-1}$ or $2^k$. 

\section{Results}

\begin{Lemma}
  \label{lemma_no_le16_4_7_without_8}
  Each $[\le 16,4,7]_2$ code contains a codeword of weight $8$.
\end{Lemma}
\begin{Proof}
  Let $C$ be an $[n,4,7]_2$ code with $n\le 16$ and $A_8=0$. From the first two MacWilliams identities we conclude
  $$
    A_7+A_9+\sum_{i\ge 10} A_i = 2^4-1=15\quad\text{and}\quad 
    7A_7+9A_9+\sum_{i\ge 10} iA_i = 2^3n =8n, 
  $$  
  so that
  $$
    2A_9+3A_{10}+\sum_{i\ge 11} (i-7)A_i = 8n-105.
  $$
  Thus, the number of even weight codewords is at most $8n/3-34$. Since at least half the codewords have to be of even weight, we obtain 
  $n\ge \left\lceil 15.75\right\rceil=16$. In the remaining case $n=16$ we use the linear programming method with the first four MacWilliams identities, 
  $A_8=0$, $B_1=0$, and the fact that there are exactly $8$ even weight codewords to conclude $A_{11}+\sum_{i\ge 13} A_i <1$, i.e., $A_{11}=0$ and $A_i=0$ for all $i\ge 13$. 
  With this and rounding to integers we obtain the bounds $5\le B_2\le 6$, which then gives the unique solution $A_7=7$, $A_9=0$, $A_{10}=6$, and $A_{12}=1$. Computing 
  the full dual weight distribution unveils $B_{15}=-2$, which is negative.
\end{Proof}

\begin{Lemma}
  \label{lemma_46_9_20_even}
  Each even $[46,9,20]_2$ code $C$ is isomorphic to a code with generator matrix
  $$
    \begin{pmatrix}
      1001010101110011011010001111001100100100000000\\
      1111100101010100100011010110011001100010000000\\
      1100110100001111101111000100000110101001000000\\
      0110101010010110101101110010100011001000100000\\
      0011101110101101100100101001010001011000010000\\
      0110011001111100011100011000110000111000001000\\
      0001111000011100000011111000001111111000000100\\
      0000000111111100000000000111111111111000000010\\
      0000000000000011111111111111111111111000000001\\
    \end{pmatrix}.
  $$
\end{Lemma}
\begin{Proof}
  Applying Proposition~\ref{prop_residual} with $w=20$ on a $[45,9,20]$ code would give a $[25,8,10]$ code, which does not exist. Thus, $C$ has effective length $n=46$, i.e., $B_1=0$. 
  Since no $[44,8,20]$ code exists, $C$ is projective, i.e., $B_2=0$. Since no $[24,8,9]$ code exists, Proposition~\ref{prop_residual} yields that $C$ cannot 
  contain a codeword of weight $w=22$. Assume for a moment that $C$ contains a codeword $c_{26}$ of weight $w=26$ and let $R$ be the corresponding residual $[20,8,7]$ 
  code. Let $c'\neq c_{26}$ be another codeword of $C$ and $w'$ and $w''$ be the weights of $c'$ and $c'+c_{26}$. Then the weight of the corresponding 
  residual codeword is given by $(w'+w''-26)/2$, so that weight $8$ is impossible in $R$ ($C$ does not contain a codeword of weight $22$). Since 
  $R$ has to contain a $[\le 16,4,7]_2$ subcode, Lemma~\ref{lemma_no_le16_4_7_without_8} shows the non-existence of $R$, so that $A_{26}=0$.

  With this, the first three MacWilliams Identities are given by
  \begin{eqnarray*}
    A_{20}+A_{24}+A_{28}+A_{30}+\sum_{i=1}^{8} A_{2i+30}&=&511\\
    3A_{20}-A_{24}-5A_{28}-7A_{30} -\sum_{i=1}^{8}\left(2i+7\right)\cdot A_{2i+30}&=&-23\\
    5A_{20}+21A_{24}-27A_{28}-75A_{30}-\sum_{i=1}^{8}\left(8i^2+56i+75\right)\cdot A_{2i+30}&=&1035.
  \end{eqnarray*}
  Minimizing $T=A_0+A_{20}+A_{24}+A_{28}+A_{32}+A_{36}+A_{40}+A_{44}$ gives $T\ge \tfrac{6712}{15}>384$, so that Proposition~\ref{prop_div_one_more}.(\ref{div_one_more_case3}) 
  gives $T=512$, i.e., all weights are divisible by $4$. A further application of the linear programming method gives that $A_{36}+A_{40}+A_{44}\le \left\lfloor\tfrac{9}{4}\right\rfloor=2$, 
  so that $C$ has to contain a $[\le 44,7,\{20,24,28,32\}]_2$ subcode.

  Next, we have used \texttt{Q-Extension} and \texttt{LinCode} to classify the $[n,k,\{20,24,28,32\}]_2$ codes for $k\le 7$ and 
  $n\le 37+k$, see Table~\ref{table_4_div_d_20}. Starting from the $337799$ doubly-even $[\le 44,7,20]$ codes, 
  \texttt{Q-Extension} and \texttt{LinCode} give 424207 doubly-even $[45,8,20]_2$ codes and no doubly-even $[44,8,20]_2$ code 
  (as the maximum minimum distance of a $[44,8]_2$ code is $19$.) Indeed, a codeword of weight $36$ or $40$ 
  can occur in a doubly-even $[45,8,20]_2$ code. We remark that largest occurring order of the automorphism group is $18$. 
  Finally, an application of \texttt{Q-Extension} and \texttt{LinCode} on the 424207 doubly-even $[45,8,20]_2$ codes results in the 
  unique code as stated. (Note that there may be also doubly-even $[45,8,20]_2$ codes with two or more codewords of a weight $w\ge 36$. 
  However, these are not relevant for our conclusion.) 
\end{Proof}

\begin{table}[htp]
  \begin{center}
    \begin{tabular}{rrrrrrrrrrrrrrrrr}
    \hline
    k / n & 20 & 24 & 28 & 30 & 32 & 34 & 35 & 36 & 37 & 38 & 39 & 40 & 41 &  42 &   43 &     44 \\ 
    \hline
    1     &  1 &  1 &  1 &  0 &  1 &  0 &  0 &  0 &  0 &  0 &    &    &    &     &      &        \\
    2     &    &    &    &  1 &  1 &  2 &  0 &  3 &  0 &  3 &  0 &    &    &     &      &        \\
    3     &    &    &    &    &    &    &  1 &  1 &  2 &  4 &  6 &  9 &    &     &      &        \\
    4     &    &    &    &    &    &    &    &    &    &  1 &  4 & 13 & 26 &     &      &        \\  
    5     &    &    &    &    &    &    &    &    &    &    &    &  3 & 15 & 163 &      &        \\
    6     &    &    &    &    &    &    &    &    &    &    &    &    &    &  24 & 3649 &        \\
    7     &    &    &    &    &    &    &    &    &    &    &    &    &    &     &    5 & 337794 \\ 
    \hline
    \end{tabular}
    \caption{Number of $[n,k,\{20,24,28,32\}]_2$ codes.}
    \label{table_4_div_d_20}
  \end{center}
\end{table}

We remark that the code of Lemma~\ref{lemma_46_9_20_even} has a trivial automorphism group and weight enumerator $1x^0+235x^{20}+171x^{24}+97x^{28}+8x^{32}$, i.e., 
all weights are divisible by four. The dual minimum distance is $3$ ($A_3^\perp=1$, $A_4^\perp=276$), i.e., the code is projective. Since the Griesmer bound, see 
Inequality~(\ref{eq_griesmer_bound}), gives a lower bound of $47$ 
for the length of a binary linear code with dimension $k=9$ and minimum distance $d\ge 21$, the code has the optimum minimum distance. The linear programming method could also be 
used to exclude the weights $w=40$ and $w=44$ directly (and to show $A_{36}\le 2$). 
While the maximum 
distance $d=20$ was proven using the Griesmer bound directly, the $[46,9,20]_2$ code is not a \textit{Griesmer code}, i.e., where Inequality~(\ref{eq_griesmer_bound}) is satisfied with equality. 
For the latter codes the $2^2$-divisibility would follow from \cite[Theorem 9]{ward2001divisible} stating that for Griesmer codes over $\mathbb{F}_p$, where $p^e$ is a divisor of the minimum 
distance, all weights are divisible by $p^e$.

\begin{Theorem}
  \label{theorem_46_9_20}
  Each $[46,9,20]_2$ code $C$ is isomorphic to a code with the generator matrix given in Lemma~\ref{lemma_46_9_20_even}.
\end{Theorem}
\begin{Proof}
  Let $C$ be a $[46,9,20]_2$ with generator matrix $G$ which is not even. Removing a column from $G$ and adding a parity check bit gives an even $[46,9,20]_2$ code. So, 
  we start from the generator matrix of Lemma~\ref{lemma_46_9_20_even} and replace a column by all $2^9-1$ possible column vectors. Checking all $46\cdot 511$ cases 
  gives either linear codes with a codeword of weight $19$ or the generator matrix of Lemma~\ref{lemma_46_9_20_even} again. 
\end{Proof}

\begin{Lemma}
  \label{lemma_no_47_10_20}
  No $[47,10,20]_2$ code exists.
\end{Lemma}
\begin{Proof}
  Assume that $C$ is a $[47,10,20]_2$ code. Since no $[46,10,20]_2$ and no $[45,9,20]_2$ code exists, we have $B_1=0$ and $B_2=0$, respectively. 
  Let $G$ be a systematic generator matrix of $C$. Since removing the $i$th unit vector and the corresponding column (with the $1$-entry) from $G$ gives a $[46,9,20]_2$ code, 
  there are at least $1023$ codewords in $C$ whose weight is divisible by $4$. Thus, Proposition~\ref{prop_div_one_more}.(\ref{div_one_more_case3}) yields that $C$ is doubly-even. 
  By Theorem~\ref{theorem_46_9_20} we have $A_{32}\ge 8$. Adding this extra inequality to the linear inequality system of the first four MacWilliams identities gives, after rounding 
  down to integers,   $A_{44}=0$, $A_{40}=0$, $A_{36}=0$, and $B_{3}=0$. (We could also conclude $B_3=0$ directly from the non-existence of a $[44,8,20]_2$-code.) The unique 
  remaining weight enumerator is given by $1x^0+418x^{20}+318x^{24}+278x^{28}+9x^{32}$. Let $C$ be such a code   and $C'$ be the code generated by the nine codewords of weight 
  $32$. We eventually add codewords from $C$ to $C'$ till $C'$ has dimension exactly nine and denote the corresponding code by $C''$. Now the existence of $C''$ contradicts 
  Theorem~\ref{theorem_46_9_20}.  
\end{Proof}

So, the unique $[46,9,20]_2$ code is strongly optimal in the sense of \cite[Definition 1]{simonis200023}, i.e., no $[n-1,k,d]_2$ and no $[n+1,k+1,d]_2$ code exists. The strongly 
optimal binary linear codes with dimension at most seven have been completely classified, except the $[56,7,26]_2$ codes, in \cite{bouyukliev2001optimal}. 
The next open case is the existence question for a $[65,9,29]_2$ code, which is equivalent to the existence of a $[66,9,30]_2$ code. The technique 
of Lemma~\ref{lemma_46_9_20_even} to conclude the $4$-divisibility of an optimal even code can also be applied in further cases and we given an example for 
$[78,9,36]_2$ codes, whose existence is unknown.    

\begin{Lemma}
  \label{lemma_no_le33_5_15_without_16}
  Each $[\le 33,5,15]_2$ code contains a codeword of weight $16$.
\end{Lemma}
\begin{Proof}
  We verify this statement computationally using \texttt{Q-Extension} and \texttt{LinCode}.
\end{Proof}
We remark that a direct proof is possible too. However, the one that we found is too involved to be presented here. Moreover, there are exactly $3$ $[\le 32,4,15]_2$ codes without 
a codeword of weight $16$.

\begin{Lemma}
  If an even $[78,9,36]_2$ code $C$ exists, then it has to be doubly-even.
\end{Lemma}
\begin{Proof}
  Since no $[77,9,36]_2$ and no $[76,8,36]_2$ code exists, we have $B_1=0$ and $B_2=0$. Proposition~\ref{prop_residual} yields that $C$ does not contain a codeword of weight $38$.  
  Assume for a moment that $C$ contains a codeword $c_{42}$ of weight $w=42$ and let $R$ be the corresponding residual $[36,8,15]_2$ code. Let $c'\neq c_{42}$ be another codeword 
  of $C$ and $w'$ and $w''$ be the weights of $c'$ and $c'+c_{42}$. Then the weight of the corresponding residual codeword is given by $(w'+w''-42)/2$, so that weight $16$ is 
  impossible in $R$ ($C$ does not contain a codeword of weight $38$). Since $R$ has to contain a $[\le 33,5,15]_2$ subcode, Lemma~\ref{lemma_no_le33_5_15_without_16} shows the 
  non-existence of $R$, so that $A_{42}=0$.

  We use the linear programming method with the first four MacWilliams identities. Minimizing the number $T$ of doubly-even codewords gives 
  $T\ge \tfrac{1976}{5}>384$, so that Proposition~\ref{prop_div_one_more}.(\ref{div_one_more_case3}) gives $T=512$, i.e., all weights are divisible by $4$.
\end{Proof}

Two cases where $8$-divisibility can be concluded for optimal even codes are given below.

\begin{Theorem}
  \label{theorem_no_85_9_40}
  No $[85,9,40]_2$ code exists.
\end{Theorem}
\begin{Proof}
  Assume that $C$ is a $[85,9,40]_2$ code. Since no $[84,9,40]_2$ and no $[83,8,40]_2$ code exists, we have $B_1=0$ and $B_2=0$, respectively. Considering the residual 
  code, Proposition~\ref{prop_residual} yields that $C$ contains no codewords with weight $w\in\{42,44,46\}$. With this, we use the first four MacWilliams identities and 
  minimize $T=A_0+\sum_{i=10}^{21} A_{4i}$. Since $T\ge 416>384$, so that Proposition~\ref{prop_div_one_more}.(\ref{div_one_more_case3}) gives $T=512$, all weights 
  are divisible by $4$. Minimizing $T=A_0+\sum_{i=5}^{10} A_{8i}$ gives $T\ge 472>384$, so that Proposition~\ref{prop_div_one_more}.(\ref{div_one_more_case3}) gives $T=512$, 
  i.e., all weights are divisible by $8$. The residual code of each codeword of weight $w$ is a projective $4$-divisible code of length $85-w$. Since no such codes of lengths 
  $5$ and $13$ exist, $C$ does not contain codewords of weight $80$ or $72$, respectively.\footnote{We remark that a $4$-divisible non-projective binary linear code of length 
  $13$ exists.} 
  
  The residual code $\hat{C}$ of a codeword of weight $64$ is a projective $4$-divisible $8$-dimensional code of length $21$. Note that $\hat{C}$ cannot contain 
  a codeword of weight $20$ since no even code of length $1$ exists. Thus we have $A_{64}\le 1$.  
  Now we look at the two-dimensional subcodes of the unique codeword of weight $64$ and two other codewords. Denoting their weights 
  by $a$, $b$, $c$ and the weight of the corresponding codeword in $\hat{C}$ by $w$ we use the notation $(a,b,c;w)$. W.l.o.g.\ we assume $a=64$,  $b\le c$ 
  and obtain the following possibilities: $(64,40,40;8)$, $(64,40,48;12)$, $(64,40,56;16)$, and $(64,48,48;16)$. Note that $(64,48,56;20)$ and $(64,56,56;24)$ are impossible. 
  By $x_8$, $x_{12}$, $x_{16}'$, and $x_{16}''$ we denote the corresponding counts. Setting $x_{16}=x_{16}'+x_{16}''$, we have that $x_i$ is the number of codewords of weight 
  $i$ in $\hat{C}$. Assuming $A_{64}=1$ the unique (theoretically) possible weight enumerator is $1x^0+360x^{40}+138x^{48}+12x^{56}+1x^{64}$. Double-counting gives 
  $A_{40}=360=2x_8+x_{12}+x_{16}'$, $A_{48}=138=x_{12}+2x_{16}''$, and $A_{56}=12=x_{16}'$. Solving this equation system gives 
  $x_{12}=348-2x_8$ and $x_{16}=x_8-93$. Using the first four MacWilliams identities for $\hat{C}$ we obtain the unique solution $x_8={102}$, $x_{12}=144$, and $x_{16}=9$, 
  so that $x_{16}''=9-12=-3$ is negative -- contradiction. Thus, $A_{64}=0$ and the unique (theoretically) possible weight enumerator is given by 
  $1x^0+361x^{40}+135x^{48}+15x^{56}$ ($B_3=60$).
  
  
  Using \texttt{Q-Extension} and \texttt{LinCode} we classify all $[n,k,\{40,48,56\}]_2$ codes for $k\le 7$ and $n\le 76+k$, see Table~\ref{table_8_div_40_48_56}. 
  For dimension $k=8$, there is no $[83,8,\{40,48,56\}]_2$ code and exactly 106322 $[84,8,\{40,48,56\}]_2$ codes. The latter codes have weight enumerators 
  $$1x^0+(186+l)x^{40}+(69-2l)x^{48}+lx^{56}$$ ($B_2=l-3$), where $3\le l\le 9$. The corresponding counts are given in Table~\ref{table_per_56}. Since the next step would need 
  a huge amount of computation time we derive some extra information on a $[84,8,\{40,48,56\}]_2$-subcode of $C$. Each of the $15$ codewords of weight $56$ of $C$ 
  hits $56$ of the columns of a generator matrix of $C$, so that there exists a column which is hit by at most $\left\lfloor 56\cdot 15/85\right\rfloor=9$ such codewords. 
  Thus, by shortening of $C$ we obtain a $[84,8,\{40,48,56\}]_2$-subcode with at least $15-9=6$ codewords of weight $56$. Extending the corresponding $5666$ cases with 
  \texttt{Q-Extension} and \texttt{LinCode} results in no $[85,9,\{40,48,56\}]_2$ code. (Each extension took between a few minutes and a few hours.)
\end{Proof}

\begin{table}[htp]
  \begin{center}
    \begin{tabular}{rrrrrrrrrrrrrrrrrrr}
    \hline
    k / n & 40 & 48 & 56 & 60 & 64 & 68 & 70 & 72 & 74 & 75 & 76 & 77 & 78 & 79 & 80 & 81 &  82 &   83 \\ 
    \hline
    1     &  1 &  1 &  1 &  0 &  0 &  0 &  0 &  0 &  0 &  0 &  0 &  0 &    &    &    &    &     &      \\
    2     &    &    &    &  1 &  1 &  2 &  0 &  2 &  0 &  0 &  2 &  0 &  0 &    &    &    &     &      \\
    3     &    &    &    &    &    &    &  1 &  1 &  2 &  0 &  3 &  0 &  5 &  0 &    &    &     &      \\ 
    4     &    &    &    &    &    &    &    &    &    &  1 &  1 &  2 &  3 &  6 & 10 &    &     &      \\ 
    5     &    &    &    &    &    &    &    &    &    &    &    &    &  1 &  3 & 11 & 16 &     &      \\
    6     &    &    &    &    &    &    &    &    &    &    &    &    &    &    &  2 &  8 & 106 &      \\
    7     &    &    &    &    &    &    &    &    &    &    &    &    &    &    &    &    &   7 & 5613 \\
    \hline
    \end{tabular}
    \caption{Number of $[n,k,\{40,48,56\}]_2$ codes.}
    \label{table_8_div_40_48_56}
  \end{center}
\end{table}

\begin{table}[htp]
  \begin{center}
    \begin{tabular}{rrrrrrrr}
    \hline
    $A_{56}$ & 3 & 4 & 5 & 6 & 7 & 8 & 9 \\ 
    \hline
     & 25773 & 48792 & 26091 & 5198 & 450 & 17 & 1\\ 
    \hline
    \end{tabular}
    \caption{Number of $[84,8,\{40,48,56\}]_2$ codes per $A_{56}$.}
    \label{table_per_56}
  \end{center}
\end{table}

\begin{Lemma}
  \label{lemma_no_le51_4_23_without_24_25_26}
  Each $[\le 47,4,23]_2$ code satisfies $A_{24}+A_{25}+A_{26}\ge 1$.
\end{Lemma}
\begin{Proof}
  We verify this statement computationally using \texttt{Q-Extension} and \texttt{LinCode}.
\end{Proof}
We remark that there a 1 $[44,3,23]_2$, 3 $[45,3,23]_2$, and 9 $[46,3,23]_2$ codes without codewords of a weight in $\{24,25,26\}$. 

\begin{Lemma}
  \label{lemma_no_even_le46_5_22_without_24}
  Each even $[\le 46,5,22]_2$ code contains a codeword of weight $24$..
\end{Lemma}
\begin{Proof}
  We verify this statement computationally using \texttt{Q-Extension} and and \texttt{LinCode}.
\end{Proof}
We remark that there a 2 $[44,4,22]_2$ and 6 $[45,4,22]_2$ codes that are even and do not contain a codeword of weight $24$.  

\begin{Theorem}
  If an even $[117,9,56]_2$ code $C$ exist, then the weights of all codewords are divisible by $8$.
\end{Theorem}   
\begin{Proof}
  From the known non-existence results we conclude $B_1=0$ and $C$ does not contain codewords with a weight in $\{58,60,62\}$. If $C$ would contain a 
  codeword of weight $66$ then its corresponding residual code $R$ is a $[51,8,23]_2$ code without codewords with a weight in $\{24,25,26\}$, which 
  contradicts Lemma~\ref{lemma_no_le51_4_23_without_24_25_26}. Thus, $A_{66}=0$. Minimizing the number $T_4$ of doubly-even codewords 
  using the first four MacWilliams identities gives $T_4\ge\frac{2916}{7}>384$, so that Proposition~\ref{prop_div_one_more}.(\ref{div_one_more_case3}) 
  gives $T_4=512$, i.e., all weights are divisible by $4$.
  
  If $C$ contains no codeword of weight $68$, then the number $T_8$ of codewords whose weight is divisible by $8$ is at least $475.86>448$, 
  so that Proposition~\ref{prop_div_one_more}.(\ref{div_one_more_case3}) gives $T_8=512$, i.e., all weights are divisible by $8$. So, let us assume 
  that $C$ contains a codeword of weight $68$ and consider the corresponding residual $[49,8,22]_2$ code $R$. Note that $R$ is even and does not contain 
  a codeword of weight $24$, which contradicts Lemma~\ref{lemma_no_even_le46_5_22_without_24}. Thus, all weights are divisible by $8$.    
\end{Proof}

\begin{Proposition}
  If an even $[118,10,56]_2$ code exist, then its weight enumerator is either $1x^0+719x^{56}+218x^{64}+85x^{72}+1x^{80}$ or $1x^0+720x^{56}+215x^{64}+88x^{72}$.
\end{Proposition}   
\begin{Proof}
  Assume that $C$ is an even $[118,10,56]_2$ code. Since no $[117,10,56]_2$ and no $[116,9,56]_2$ code exists we have $B_1=0$ and $B_2=0$, respectively. Using 
  the known upper bounds on the minimum distance for $9$-dimensional codes we can conclude that no codeword as a weight $w\in\{58,60,62,66,68,70\}$. Maximizing 
  $T=\sum_i A_{4i}$ gives $T\ge 1011.2>768$, so that $C$ is $4$-divisible, see Proposition~\ref{prop_div_one_more}.(\ref{div_one_more_case3}). Maximizing $T=\sum_i A_{8i}$ 
  gives $T\ge 1019.2>768$, so that $C$ is $8$-divisible, Proposition~\ref{prop_div_one_more}.(\ref{div_one_more_case3}). 
  Maximizing $A_{i}$ for $i\in \{88,96,104,112\}$ gives a value strictly less than $1$, so that the only non-zero weights can be $56$, $64$, $72$, and $80$. 
  Maximizing $A_{80}$ gives an upper bound of $\frac{3}{2}$, so that $A_{80}=1$ or $A_{80}=0$. The remaining values are then uniquely determined by the first four 
  MacWilliams identities. 
\end{Proof}

The exhaustive enumeration of all $[117,9,\{56,64,72\}]_2$ codes remains a computational challenge. While we have constructed a few thousand non-isomorphic 
$[115,7,\{56,64,72\}]_2$ codes, we still do not know whether a $[117,9,56]_2$ code exists.


\end{document}